\documentclass{article}
\usepackage{amsfonts,amsmath,graphicx,color}
\usepackage{psfrag}

\newcommand{\R}{\mathbb{R}}
\newcommand{\C}{\mathbb{C}}

\begin{document}

\title{Tipping points in open systems: bifurcation, noise-induced and rate-dependent examples in the climate system}

\author{Peter Ashwin, Sebastian Wieczorek, Renato Vitolo \& Peter Cox\\
Mathematics Research Institute,\\
University of Exeter,\\
Exeter EX4 4QF.
}

\maketitle

\begin{abstract}
  Tipping points associated with bifurcations (B-tipping) or induced by noise
  (N-tipping) are recognized mechanisms that may potentially lead to sudden
  climate change. We focus here a novel class of tipping points, where 
  a sufficiently rapid change to an input or parameter of a system may cause the
  system to ``tip'' or move away from a branch of attractors. Such rate-dependent
  tipping, or {\em R-tipping}, need not be associated with either bifurcations 
  or noise.
  We present an example of all three types of tipping in a  
  simple global energy balance model of the climate system, illustrating the 
  possibility of dangerous rates of change even in the absence of noise and of 
  bifurcations in the underlying quasi-static system.
\end{abstract}

~

\section{Tipping points - not just bifurcations}

In the last few years, the idea of ``tipping points'' has caught the
imagination in climate science with the possibility, also indicated by both
palaeoclimate data and global climate models, that the climate system may
abruptly ``tip'' from one regime to another in a comparatively short time.\footnote{NB This version of the paper includes a correction to Section 2.1.}

This recent interest in tipping points is related to a long-standing question
in climate science: to understand whether climate fluctuations and transitions
between different ``states'' are due to external causes (such as variations in
the insolation or orbital parameters of the Earth) or to internal mechanisms
(such as oceanic and atmospheric feedbacks acting on different timescales).  A
famous example is Milankovich theory, according to which these transitions are
forced by an external cause, namely the periodic variations in the Earth's
orbital parameters (see e.g.~\cite{Hays1976}). Remarkably, the evidence in
favour of Milankovich theory still remains controversial, see
e.g.~\cite{Shackleton2000}.

Hasselmann \cite{Hasselmann1976} was one of the first to tackle this question
through simple climate models obtained as stochastically perturbed dynamical
systems. He argued that the climate system can be conceptually divided into
a fast component (the ``weather'', essentially corresponding to the evolution
of the atmosphere) and a slow component (the ``climate'', that is the ocean,
cryosphere, land vegetation, etc.). The ``weather'' would act as an
essentially random forcing exciting the response of the slow ``climate''.  In
this way, short-time scale phenomena, modelled as stochastic perturbations,
can be thought of as driving long-term climate variations.  This is what we
refer to as \emph{noise-induced tipping}.

Sutera~\cite{Sutera1981} studies noise-induced tipping in a simple global energy balance
model previously derived by Fraedrich~\cite{Fraedrich1979}.  Sutera's results
indicate a characteristic time of $10^5$yr for the the random transitions
between the ``warm'' and the ``cold'' climate states, which matches well with
the observed average value. One shortcoming is that this analysis leaves open
the question as to the \emph{periodicity} of such transitions indicated by the
power spectral analysis \cite[Fig. 7]{Mason1976}.

There is a considerable literature on noise-induced escape from attractors in stochastic models \cite{Hanggi}. These have successfully been used for
modelling changes in climate phenomena \cite{Sura2002}, although authors
do not always use the word ``tipping'' and other aspects have been examined. For example, Kondepudi {\em et al} \cite{Kondepudi_etal_1986} consider the combined effect of noise and parameter changes on the related problem of ``attractor selection'' in a noisy system.

More recently, bifurcation-driven tipping points or {\em dynamic bifurcations} \cite{DynBifs} have been suggested as an important mechanism by which sudden changes in the behaviour of a system may come about. For example, Lenton {\em et al}~\cite{Lenton_etal_2008,Lenthisvol}
conceptualize this as an open system
\begin{equation}
  \label{eq:xdot} 
  \frac{dx}{dt}=f(x,\lambda(t))
\end{equation}
where $\lambda(t)$ is in general a time-varying input. 
In the case that $\lambda$ is constant, we refer to (\ref{eq:xdot}) as the {\em parametrized system} 
with parameter $\lambda$, and to its stable solution as the {\it quasi-static attractor}. 
If $\lambda(t)$ passes through a bifurcation point of the parametrized system where a quasi-static 
attractor (such as an equilibrium point $\tilde{x}(\lambda)$) loses stability, it is
intuitively clear that a system may ``tip'' directly as a result of varying that
parameter, though in certain circumstances the effect may be delayed because of well-documented slow passage through bifurcation effects \cite{BifDelay}. Related to this, Dakos et al \cite{Dakos_etal_2008}
have proposed that tipping points are recognizable and to some extent
predictable. They propose a method to de-trend signals and then,
examining the correlation of fluctuations in the detrended signal, they find a signature of bifurcation-induced tipping points. These papers have concentrated on systems where equilibrium solutions for the parametrised system lose stability, although recent work of Kuehn \cite{Kuehn} considers tipping effects in general two timescale systems as occasions where there is a bifurcation of the fast dynamics. 

The explanation of climate tipping as a phenomenon purely induced by
bifurcations has been called into question. For example, Ditlevsen and
Johnsen~\cite{Ditlevsen_etal_2010} suggest that the predictive techniques
to forecast a forthcoming tipping point~\cite{Dakos_etal_2008} are not always
reliable. Indeed, noise alone can drive a system to tipping without any
bifurcation. Nonetheless, it seems that the ideas of bifurcation-induced
tipping can give practically useful predictions, for example in detecting
potential ecosystem population collapses \cite{Drake_etal_2010}.

In their review paper, Thompson and Sieber~\cite{Thompson_Sieber_2010b}
discuss bifurcation- and noise-induced mechanisms for tipping.  They examine
stochastically forced systems
\begin{equation} \label{eq:xdotstoch}
dx=f(x,\lambda(t))dt+g(x)dW
\end{equation}
where $W$ represents a Brownian motion. Using generic bifurcation theory, they
distinguish between \emph{safe bifurcations} (where an attracting state loses
stability but is replaced by another ``nearby'' attractor), \emph{explosive
  bifurcations} (where the attractor dynamics explores more of phase space but
still returns to near the old attractor) and \emph{dangerous bifurcations}
(where the attractor dyamics after bifurcation are unrelated to what has gone
before).  In~\cite{Thompson_Sieber_2010a}, Thompson and Sieber clarify that a
timeseries analysis of bifurcation-induced tipping point near a quasi-static
equilibrium (QSE) relies on a separation of timescales
\begin{equation}
\kappa_{\mathrm{drift}}\ll\kappa_{\mathrm{crit}}\ll\kappa_{\mathrm{stab}}
\label{eq:rates}
\end{equation}
where $\kappa_{\mathrm{drift}}$ is the average drift rate of parameters,
$\kappa_{\mathrm{crit}}$ is the decay rate for the slowest decaying mode of
the QSE and $\kappa_{\mathrm{stab}}$ are the remaining (faster decaying)
modes. However, it is not easy to define $\kappa_{\mathrm{drift}}$ in general
(especially in a coordinate-independent manner) and there is no {\em a priori}
reason for inequality (\ref{eq:rates}) to hold for a given system.  

Rate-dependent tipping has not previously been discussed in detail, but it has been identified in \cite{Wieczorek_etal_2010}  as an important tipping mechanism that cannot be explained by previously proposed mechanisms ({\it i.e.} noise or bifurcations of a quasi-static attractor). This paper aims to better understand the phenomenon of rate-dependent tipping by introducing a linear model with a {\it tipping radius} and discussing three basic examples where this type of tipping appears. 

We suggest that tipping effects in open systems can be usefully split into three categories:
\begin{itemize}
\item ``B-tipping'' where the output from an open system changes abruptly or
  qualitatively due to a bifurcation of  a quasi-static attractor.
\item ``N-tipping'' where noisy fluctuations result in the system departing
  from a neighbourhood of a  quasi-static attractor.
\item ``R-tipping''  where the system fails to track a continuously changing quasi-static attractor.
\end{itemize}
 We demonstrate that each mechanism {\em on its own} can produce a tipping response. Furthermore, any open system may exhibit tipping phenomena that result from a
combination of several of the above. 

The paper is organized as follows; in the remainder of this section we discuss
a setting for open systems allowing discussion of the three types of tipping
phenomena. In Section~\ref{sec:R-tipping} we formulate a criterion for
R-tipping. In Section~\ref{sec:R-tipping-examples} we discuss three illustrative low dimensional examples of R-tipping; two related to bifurcation normal forms and one for a slow-fast system. Section~\ref{sec:BNR-tipping_egs} gives an illustrative examples of all three types of tipping for an energy-balance model of the global climate for different parameter regimes. Section~\ref{sec:discuss} concludes with a discussion and some open questions.


\subsection{Towards a general theory of tipping in open systems}

Dynamical systems theory has developed a wide-ranging corpus of results concentrated on the behaviour of autonomous finite dimensional deterministic systems - often called {\em closed systems}, because their future time trajectories depend only on the current state of the system. If the systems have inputs that can change the fate of system trajectories then we say the system is {\em open}. Real-world systems are never closed except to some degree of approximation, and a range of methods have been developed to cope with the fact that they are open:
(a) One can view external perturbations as time variation of parameters that would be fixed for a closed system. (b) There are various theoretical approaches to stochasticity in systems, either intrinsic or external. (c) Control theory allows one to design inputs to control a system's outputs in a desired way, given (possibly imperfect) knowledge of the system.

\begin{figure}%
\begin{center}
\includegraphics[width=6cm,clip=]{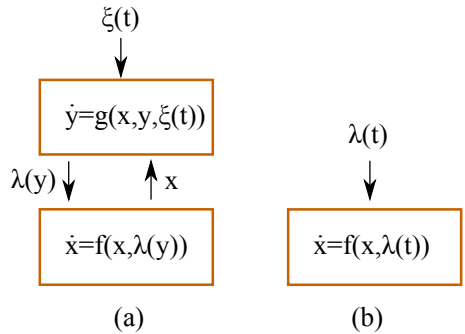}%
\end{center}
\caption{The system $(x,y)$ in (a) represents a (high dimensional) open system with inputs $\xi(t)$. We identify a low dimensional but nonlinear subsystem $x$ forced by some ``observables'' $\lambda(y)$ from the high dimensional system. The behaviour of $x$ in (a) can be partially understood by examining the open subsystem (b) for a suitable class of temporal forcing $\lambda(t)$.}%
\label{fig:blockdiagram}%
\end{figure}

\begin{figure}%
\begin{center}
\includegraphics[width=13cm,clip=]{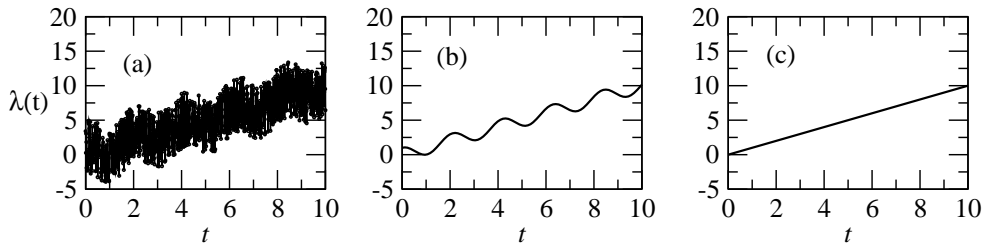}%
\end{center}
\caption{Candidate inputs $\lambda(t)$ in the subsystem {\protect (\ref{eq:xdot})}. These may include noise- and signal-like components as in (a), or purely deterministic smooth inputs/parameter variations as in (b,c). Tipping responses in the subsystem may occur in response to noise fluctuations (N-tipping), to passing through a bifurcation point for the parametrized subsystem $\lambda=\lambda_0$ (B-tipping) or as a result of too rapid variation (R-tipping). }%
\label{fig:timeseries}
\end{figure}

In Fig.~\ref{fig:blockdiagram}(a) we illustrate an arbitrary high dimensional system  where we have identified a low dimensional subsystem that we wish to check for ``tipping effects''. We do this by analysing the response of an open system (\ref{eq:xdot}) in Fig.~\ref{fig:blockdiagram}(b) to possible time-varying inputs $\lambda(t)$. Fig.~\ref{fig:timeseries} shows some possible candidates for the input $\lambda(t)$; we are interested in classifying those inputs that lead to a sudden change in $x$. This ``tipping'' may depend on details of the noise (N-tipping), may involve passing through a critical value of $\lambda(t)$ corresponding to a bifurcation (B-tipping) of the parametrized subsystem or may depend on the rate of change of $\lambda(t)$ along some path in parameter space (R-tipping).

\section{R-tipping: a linear model}
\label{sec:R-tipping}

%

We use a simple model to explore R-tipping and to give sufficient conditions such that R-tipping does/does not occur. Suppose that the system (\ref{eq:xdot}) for $x\in\R^n$ and parameter $\lambda$ has a {\em quasi-static equilibrium} (QSE) $\tilde{x}(\lambda)$ with a {\em tipping radius} $R>0$. For some initial $x_0$ with $|x_0-\tilde{x}(\lambda)|<R$ we assume that the evolution of $x$ with time is given by
\begin{equation}\label{eq:model}
\frac{dx}{dt}=M(x-\tilde{x}(\lambda))~~\mbox{ for }~~|x-\tilde{x}(\lambda)|<R
\end{equation}
where $M$ is a fixed stable linear operator (i.e. $|e^{Mt}|\rightarrow 0$ as $t\rightarrow \infty$).
More generally, we consider a time varying parameter, $\lambda(t)$, that represents the input to the subsystem. 
If $|x(t)-\tilde{x}(\lambda(t))|<R$ then we say that $x(t)$ {\em tracks} (or adiabatically follows) the QSE $\tilde{x}(\lambda)$. If there is a $t_0$ such that $|x(t_0)-\tilde{x}(\lambda(t_0))|= R$ then we say the solution {\em tips} (adiabatic approximation fails) at $t_0$ and regard the model as unphysical beyond this point in time.
The tipping radius may be related to the basin of attraction boundary for the nonlinear problem~(\ref{eq:xdot}), as is the case in Sec.~3(a,b), but it need not be, as is the case in Sec.~3(c) and~\cite{Wieczorek_etal_2010}. System (\ref{eq:model}) shows only R-tipping -  because $M$ is fixed there is no bifurcation in the system and no noise is present. Clearly the model can be generalised to include $M$ and $R$ that vary with $\lambda(t)$, and/or nonlinear terms. 
 Equation (\ref{eq:model}) can be solved with initial condition $x(0)=x_0$ to give
$$
x(t)= e^{Mt} x_0 - \int_{s=0}^t e^{M(t-s)} M \tilde{x}(\lambda(s))\, ds.
$$
If we assume that the solution is modelled by the linear system near the QSE for an arbitrarily long past and set $u=t-s$ then the dependence on initial value decays to give
\begin{equation}
x(t)= -\int_{u=0}^{\infty} e^{Mu} M \tilde{x}(\lambda(t-u))\, du.
\label{eq:xt1}
\end{equation}
Assuming that $M$ is invertible and exponentially stable (more precisely, we assume that $|e^{Mt}M^{-k}v|\rightarrow 0$ as $t\rightarrow 0$ for $k=1,2$) and that the rate of motion of the QSE and parameter are bounded (more precisely, the derivatives $d^l\tilde{x}/d\lambda^l$ and $d^l\lambda/dt^l$ for $l=1,2$ are bounded in time) then (\ref{eq:xt1}) can be integrated by parts to give
$$
x(t) = -\left[e^{Mu} \tilde{x}(\lambda(t-u))\right]_0^{\infty}
-\int_0^{\infty} e^{Mu} \frac{d\tilde{x}}{dt}(\lambda(t-u)) \,du
$$
and so
\begin{equation}\label{eq:xminustildex}
x(t) -\tilde{x}(\lambda(t)) = - \int_0^{\infty} e^{Mu} \frac{d\tilde{x}}{dt}(\lambda(t-u))\,du.
\end{equation}
Integrating again by parts gives
\begin{eqnarray*}
x(t) -\tilde{x}(t) &=& M^{-1}\frac{d\tilde{x}}{dt}(\lambda(t-u))
 -\int_0^{\infty} e^{Mu}M^{-1} \frac{d^2\tilde{x}}{dt^2}(\lambda(t-u))\,du\\
&=& L(t)+E(t)
\end{eqnarray*}
The linear instantaneous lag is
\begin{equation}
L(t)= M^{-1} \frac{d\tilde{x}}{dt}(\lambda(t))
\end{equation}
If we define the {\em drift} of the QSE to be the rate of change
\begin{equation}\label{eq:rdrift}
r(t):=\frac{d\tilde{x}}{dt}=\frac{d\tilde{x}}{d\lambda}\,\frac{d\lambda}{dt}
\end{equation}
then the linear instantaneous lag is
\begin{equation}
L(t)= M^{-1} r(t).
\label{eq:linlag}
\end{equation}
The error to the linear instantaneous lag is
$$
E(t)=-\int_0^{\infty} e^{Mu}M^{-1} \frac{d^2\tilde{x}}{dt^2}(\lambda(t-u)) \,du.
$$
which includes historical information. This can also be expressed as
$$
E(t)=-\int_0^{\infty} e^{Mu}M^{-1} \left[\tilde{x}''(\lambda(t-u))(\lambda'(t-u))^2+\tilde{x}'(\lambda(t-u))\lambda''(t-u)\right]\,du.
$$
To summarize, the solution of (\ref{eq:model}) follows the QSE $\tilde{x}(\lambda(t))$ with a linear instantaneous lag term $L(t)$ and a history-dependent term $E(t)$.

\subsection{A criterion for R-tipping with steady drift}

If $\frac{d\tilde{x}}{dt}=r$ is constant in time then we say the system has {\em steady drift} and (\ref{eq:xminustildex}) simplifies to $x(t)-\tilde{x}(\lambda(t))=M^{-1}r$. In other words, one can verify that $E(t)=0$ and that
\begin{equation}
|x(t)-\tilde{x}(\lambda(t))|=|M^{-1}r|,
\label{eq:dev_from_qe}
\end{equation}
On writing the matrix norm $\|M\|=\sup_{v\neq 0} |Mv|/|v|$ we note that for any $r\neq 0$ and invertible $M$ we have
$$
\|M\|^{-1}\cdot|r|\leq |M^{-1}r|\leq \|M^{-1}\|\cdot|r|.
$$
We can avoid R-tipping if $|x(t)-\tilde{x}(\lambda(t))|=|M^{-1}r|< R$ and hence a sufficient condition on the rate of parameter variation to avoid R-tipping is that
\begin{equation}
 \|M^{-1}\| \cdot |r| < R
\label{eq:r_tip_drift}
\end{equation}
while a sufficient condition for R-tipping to occur in this model is that
$$
\|M\|^{-1}\cdot |r|>R.
$$
In the intermediate case, the path of parameter variation will determine whether or not there is any R-tipping.

\subsection{General criteria for R-tipping}

In the more general case\footnote{The first part of this section is a revised and corrected version of the corresponding section that appeared in in Phil. Trans. R. Soc. A (2012) 370, 1166-1184, and is based on the published correction prepared with the assistance of Clare Hobbs.} where $r(t)$ varies with $t$ we use (\ref{eq:xminustildex}) and the inequality $|e^{Mu} v|\leq \|e^{Mu}\| |v|$ to give the upper bound
$$
|x(t)-\tilde{x}(t)|\le \int_0^{\infty}\left\| e^{Mu}\right\|\left|\frac{d\tilde{x}}{dt}(\lambda(t-u))\right| \,du. 
$$
If we define
$$
r_{max}(t)=\sup_{s\leq t} \left|\frac{d\tilde{x}}{dt}(\lambda(s))\right|
           =\sup_{s\leq t} \left|\frac{d\tilde{x}}{d\lambda}(\lambda(s))\;\frac{d\lambda}{dt}(s)\right|,
$$
then this means that
$$
|x(t)-\tilde{x}(t)|\le r_{max}(t)\int_0^{\infty}\left\| e^{Mu}\right\|\,du.
$$
Moreover, if $M$ is stable then
\begin{eqnarray}
\left\| e^{Mu}\right\|\le c\,e^{-\beta u}
\label{eq:upperbound}
\end{eqnarray}
for some real $c,\beta>0$  (see \cite{Hinrichsen2011}), and so
$$
|x(t)-\tilde{x}(t)|\le r_{max}(t)\, \frac{c}{\beta}.
$$ 
Hence, we can guarantee that (\ref{eq:model}) avoids R-tipping by time $t$ if
\begin{eqnarray}
 \frac{c}{\beta}\,r_{max}(t) <R.
\label{eq:avoids}
\end{eqnarray}
If $M$ is scalar then we can choose $c=1$, $\beta=-M$ and (\ref{eq:avoids}) reduces to
\begin{equation}
\|M^{-1}\| \cdot r_{max}(t) <R
\label{eq:r_tip_drift_u}
\end{equation}
On the other hand, if $M$ is a matrix then we need a good choice of $c$ and $\beta$ in~(\ref{eq:upperbound}) to make the tipping condition~(\ref{eq:avoids}) optimal, but this depends on the matrix structure and not simply the norm; see for example the text by Hinrichsen and Pritchard \cite{Hinrichsen2011} and the elegant estimates of Godunov~\cite[Eq.(13)]{Veselic2003}. 

Note that the tipping radius approach does not seem to easily give any sufficient condition to undergo R-tipping in terms of $r_{max}(t)$, comparable to sufficient condition for tipping with steady drift given the end of the previous section.

One can define a natural timescale for the motion of the QSE as
$$
\frac{R}{|r(t)|}=R\left(\left|\frac{d\tilde{x}}{d\lambda}\frac{d\lambda}{dt}\right|\right)^{-1};
$$
note that  in general, combinations of $d\tilde{x}/dt$ and $d\lambda/dt$ do not give timescales in units s$^{-1}$. For an R-tipping to occur, this natural timescale 
may be comparable to the slowest timescale (e.g. the reciprocal of the leading eigenvalue of $M$)  of the parametrized system. 
The three examples in the next section have $ |d\tilde{x}/d\lambda|=1$ and $R\approx 1$ so we expect 
R-tipping when $|d\lambda/dt|\approx\|M\|$. However, if $ |d\tilde{x}/d\lambda|\approx 1/\epsilon$, then 
clearly we can have R-tipping even when $|d\lambda/dt|\approx\epsilon\|M\|$. 

It is possible to think of more general tipping problems by analogy with the ``linear system and tipping radius'' model discussed here. For example, for an open nonlinear system we consider an ``effective tipping radius'' that corresponds to how far the linearized system needs to be from a branch of QSE to ensure that the nonlinear system tips. There is however no exact analogy possible - the effective tipping radius may depend on the shape of the local basin of attraction, the nonlinearities present and the exact path taken in parameter space.


\section{R-tipping: model examples}
\label{sec:R-tipping-examples}

We give three illustrative examples of R-tipping. The first two are based on normal forms for the two basic codimension one bifurcation that are generic for dissipative systems: the saddle-node and the Hopf bifurcation. The third is an example that uses a fast-slow system to show that R-tipping can occur even in cases where there is a single attractor that is globally asymptotically stable. 

\subsection{Saddle-node normal form with steady drift}
\label{sec:sndrift}

\begin{figure}
  \begin{center} \includegraphics[type=eps,ext=.eps,read=.eps,width=13cm]
    {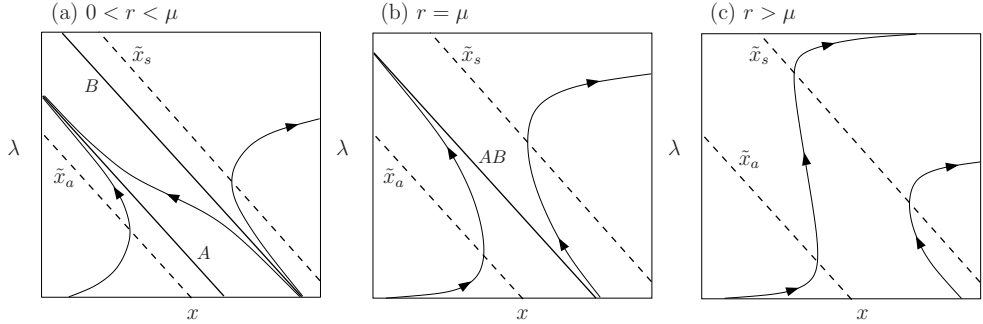}
    \caption{\label{fig:A}
Phase portraits of~(\ref{eq:nf2}) for (a) $0<r<\mu$, (b) $r=\mu$, and (c) $r>\mu$, including the two quasi-static equilibria, $\tilde{x}_a$ and $\tilde{x}_s$, and the two invariant lines, $A$ and $B$.
}
\end{center}
\end{figure}

We consider the example system for $x\in\R$ with parameter $\lambda(t)\in \R$ and drift $r$.
\begin{eqnarray}
\label{eq:nf2}
\frac{dx}{dt} &=& \left(x+\lambda\right)^2 -\mu,\\
\label{eq:nf2r}
\frac{d\lambda}{dt} &=& r,
\end{eqnarray}
with fixed $\mu>0$. In the $(x,\lambda)$ phase plane of~(\ref{eq:nf2})--(\ref{eq:nf2r}), there are two $dx/dt=0$ isoclines given by
$
\tilde{x}_a(\mu)=\{(x,\lambda)\in \mathbb{R}^{2}:\lambda=-\sqrt{\mu}-x\}
$
and
$
\tilde{x}_s(\mu)=\{(x,\lambda)\in \mathbb{R}^{2}:\lambda=\sqrt{\mu}-x\},
$
that correspond to two QSE; a stable node and a saddle, respectively, for~(\ref{eq:nf2}).
Furthermore, if $\mu>r$ there are two invariant lines, one attracting
$$
A(\mu,r)=\{(x,\lambda)\in \mathbb{R}^{2}:\lambda=-\sqrt{\mu-r}-x\},
$$
and one repelling 
$$
B(\mu,r)=\{(x,\lambda)\in \mathbb{R}^{2}:\lambda=\sqrt{\mu-r}-x\},
$$
both with a constant slope $d\lambda/dx =-1$ (Fig.~\ref{fig:A}). The stability manifests itself as an exponential decay (growth) of a small perturbations about  $A(\mu,r)$ $(B(\mu,r))$.

If $0<r<\mu$ then $B(\mu,r)$ defines a {\em tipping threshold}. Initial conditions below $B(\mu,r)$ converge to $A(\mu,r)$ whereas initial conditions above $B(\mu,r)$ give rise to solutions $x(t)\rightarrow\infty$ as $t\rightarrow\infty$. If $r=\mu$ then $A(\mu,r)$ and $B(\mu,r)$ coalesce into a neutrally stable invariant line $AB$ [Fig.~\ref{fig:A}(b)] that disappears  for $r>\mu$ [Fig.~\ref{fig:A}(c)]. Hence, for $r>\mu$ there is no tipping threshold, meaning that trajectories for all initial conditions become unbounded as $t\rightarrow\infty$.

Let us assume that the system is at $(x_0,\lambda_0)$ at time $t=0$. If the initial condition $(x_0,\lambda_0)$ lies between $\lambda=-x$ and $\tilde{x}_s(\mu)$ then the critical rate $r_c$ is the value of $r$ at which  the $r$-dependent  tipping threshold $B(\mu,r)$ crosses $(x_0,\lambda_0)$. If the initial condition lies on or below the line $\lambda=-x$ then  the critical rate $r_c$ is the value of $r$ at which $B(\mu,r)$ and $A(\mu,r)$ meet and disappear. This gives a precise value for the critical rate as the following function of initial conditions:
\begin{eqnarray}
\label{eq:crit}
r_c=\left\{
\begin{array}{ll}
\mu-(\lambda_0+x_0)^2 &~~\mbox{if}~~~-x_0<\lambda_0<-x_0+\sqrt{\mu},\\
\mu &~~\mbox{if}~~~~ \lambda_0\le -x_0.
\end{array}\right.
\end{eqnarray}

We can approximate this result using the simple linear model (\ref{eq:model}) with the linearization at the QSE as $M=-2\sqrt{\mu}$ so that $\|M^{-1}\|=\|M\|^{-1}=1/(2\sqrt{\mu})$. Clearly, the linear model with $R=2\sqrt{\mu}$ given by $\tilde{x}_s$ (basin boundary for $\tilde{x}_a$) overestimates $r_c$ because  the linear attraction weakens on moving away from the stable QSE in the nonlinear problem. This can be overcome by choosing an {\em effective tipping radius} $R_c$. Comparing with (\ref{eq:r_tip_drift}), the system avoids tipping if
$$
|r|< 2\sqrt{\mu} R_c
$$
which, for $r_c=\mu$, suggests an effective tipping radius $R_c=\sqrt{\mu}/2$. 
Finally, owing to steady drift, this problem can  be reduced to a saddle-node bifurcation at $r=\mu$ in a co-moving coordinate system $y=x+\lambda$.

\subsection{A subcritical Hopf normal form}

As a second example  we consider 
\begin{equation}
\frac{dz}{dt}=F(z-\lambda),
\label{eq:hopf}
\end{equation}
where $z=x+iy\in \C$.  For the subcritical Hopf normal form with frequency $\omega$ we choose
$$
F(z)=(-1+ i\,\omega)z+|z|^2z.
$$ 
Note that the system (\ref{eq:hopf}) has only one QSE at $\tilde{z}=\lambda(t)$. Two cases of R-tipping that we consider are with steady drift ( these can be reduced to a bifurcation problem in another coordinate system) and with unsteady drift (where there is no straightforward simplification to a bifurcation problem).

\subsubsection{Hopf normal form with steady drift}

\begin{figure}%
\begin{center}
\includegraphics[height=6.3cm,clip=]{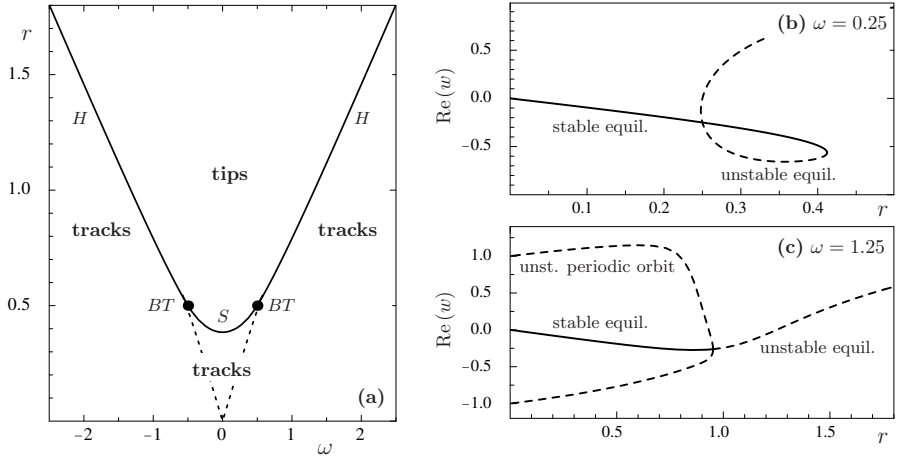}%
\end{center}
\caption{
(a) Solid curves in a two-parameter tipping diagram for (\ref{eq:hopf}) with steady drift
indicate the critical rate $r_c(\omega)$. The stable equilibrium for the co-moving system (\ref{eq:hopfrate}), or the ability to track the QSE in the original system (\ref{eq:hopf}), (b) disappears in a saddle-node bifurcation  or (c) destabilises in a subcritical Hopf bifurcation when $r=r_c(\omega)$.
}
\label{fig:td}
\end{figure}

Consider (\ref{eq:hopf}) with a uniform drift of the QSE along the real axis, at a rate $r$ (which must be real): 
$
d\lambda/dt=r.
$
There is a critical rate $r_c$ at which the system moves away from the stable QSE. We can find this $r_c$ analytically by changing to the co-moving system for $w=z-\lambda$,
\begin{equation}
\frac{dw}{dt}=F(w)- r,
\label{eq:hopfrate}
\end{equation} 
where a stable equilibrium represents the ability to track the QSE in the original system.
Setting $w=|w|\,e^{i\theta}$ and  rewriting equation (\ref{eq:hopfrate}) in terms of  $d|w|/dt$ and $d\theta/dt$ gives an equilibrium at $(|w_e|,\theta_e)$ that satisfies
\begin{equation}
|w_e|^6-2|w_e|^4+(\omega^2+1)|w_e|^2-r^2=0.
\label{eq:equilcond2}
\end{equation}
In the $(r,\omega)$ parameter plane, there is a saddle-node bifurcation curve ($S$ in  Fig.~\ref{fig:td}(a))
whose different branches are given by equation (\ref{eq:equilcond2}) with
\begin{equation}
|w_e|^2_\pm=\frac{2}{3}\left(1\pm\sqrt{1-\frac{3}{4}(1+\omega^2)} \right),
\label{eq:sn}
\end{equation}
and join at cusp points at $(r,\omega)=(\pm(2/3)^{3/2},\pm(1/3)^{1/2})$ 
(not marked in Fig.~\ref{fig:td}(a)).
Linearising about the stable equilibrium $(|w_e|,\theta_e)$ reveals that the characteristic 
polynomial
$$
s^2+(2-4|w_e|^2)\,s + \omega^2 +(|w_e|^2-1)(3|w_e|^2-1)=0,
$$
has a pair of pure imaginary roots indicating a Hopf bifurcation when
$|w_e|^2=1/2$ and $\omega^2>1/4$.
In the $(r,\omega)$ parameter plane, (disjoint) Hopf bifurcation curves ($H$) originate 
from Bogdanov-Takens bifurcation points ($BT$) at $(r,\omega)=(\pm 1/2,\pm 1/2)$,  and 
are given by
\begin{equation}
\frac{1+4\omega^2}{8}-r^2=0\;\;\;\;\mbox{and}\;\;\;\;\omega^2>1/4,
\label{eq:Hopf}
\end{equation} 
that follows from equation (\ref{eq:equilcond2}) with $|w_e|^2=1/2$.
At $BT$, saddle-node bifurcation changes from super (solid) to subcritical (dashed).
It turns out that the stable equilibrium for (\ref{eq:hopfrate}), indicating the  
ability to track the QSE in the original system,  disappears in a supercritical saddle-node 
bifurcation when $\omega^2<1/4$ and becomes unstable in a subcritical Hopf bifurcation 
when $\omega^2>1/4$. Hence, for initial conditions within the basin boundary of this 
equilibrium, the critical rate is given by
\begin{eqnarray}
\label{eq:critHopf}
r_c(\omega)=\left\{
\begin{array}{ll}
\pm\sqrt{|w_e|_-^6-2|w_e|_-^4+(\omega^2+1)|w_e|_-^2} &~~\mbox{if}~~~\omega^2\le 1/4,\\
\pm\sqrt{(1+4\omega^2)/8} &~~\mbox{if}~~~\omega^2>1/4.
\end{array}\right.
\end{eqnarray} 

Again, we can approximate this result using the simple linear model (\ref{eq:model}) with the linearization at the QSE as $M=(-1,5;-5,-1)$ so that $\|M^{-1}\|=\|M\|^{-1}=0.1961$. Clearly, the linear model with a tipping radius $R=1$ given by the unstable periodic orbit (basin boundary for $\tilde{z}$) does not account for nonlinear attraction away from the QSE and for the spiraling shape of trajectories when $\omega\ne0$. Therefore, we choose an effective tipping radius $R_c$. Comparing with (\ref{eq:r_tip_drift}), the system avoids tipping if
$$
|r|< 5.0990 R_c
$$
which suggests an $\omega$-dependent effective tipping radius $R_c(\omega)=r_c(\omega)/5.099$. 

R-tipping that reduces to a bifurcation problem in a co-moving system should not 
be confused with B-tipping: observe that the bifurcation parameter $r$ does not 
vary in time,  and it is ``the ability to track the QSE'', rather than the QSE 
itself, that bifurcates.

\subsubsection{Hopf normal form with unsteady drift}

\begin{figure}
\begin{center}
\includegraphics[width=10cm,clip=]{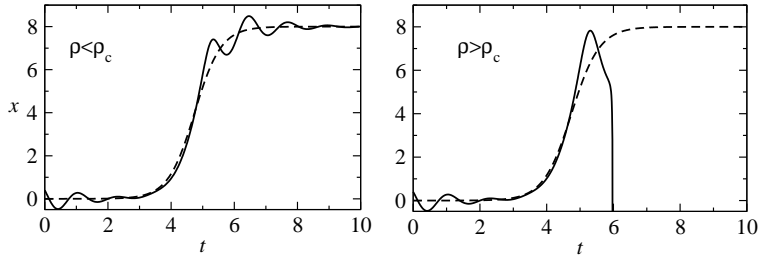}%
\end{center}
\caption{R-tipping for ({\protect \ref{eq:hopfrate})--(\ref{eq:shift}}) for $\Delta=8$ showing time evolution for (a) $\rho=4.76$ and (b) $\rho=4.8$ (recall that $\rho$ scales the maximum rate of change) from an initial condition $(x,y,\lambda)=(0.4,0.5,0.0001)$. For $\rho>\rho_c=4.78$ we find that system trajectories no longer follow the stable quasi-static equilibrium (shown by the dashed line) as they meet its basin boundary.
}%
\label{fig:hopf_shift_rate}%
\end{figure}

We now consider (\ref{eq:hopf}) where we include a smooth shift of QSE between asymptotically steady positions at $z=0$ to $z=\Delta$, according to
\begin{equation}
\frac{d\lambda}{dt}= \rho \lambda(\Delta-\lambda),~~\lambda(t_0)=\Delta/2
\label{eq:shift}
\end{equation}
where $\rho>0$ parametrizes the {\em maximum rate of the shift}, $\Delta>0$ is the amplitude of the shift and $t_0$ is the time when the rate of change is largest. Integrating (\ref{eq:shift}) gives
\begin{equation}
\lambda(t)= \Delta(\tanh(\Delta\;\rho (t-t_0)/2)+1)/2.
\label{eq:shiftpath}
\end{equation}
which implies the following parameter dependence on time:
$$
\lambda(-t)\rightarrow 0,~~\lambda(t)\rightarrow \Delta~~\mbox{ as }~t\rightarrow \infty~\mbox{ and }~
\frac{d\lambda}{dt}\leq\frac{d\lambda}{dt}(t_0)=\frac{\Delta^2\,\rho}{4}.
$$
Near $t=t_0$ this describes a smooth shift between the location of an asymptotically stable equilibrium from $z=0$ to $z=\Delta$, and the maximum rate of the shift is $\Delta^2\,\rho/4$ at $t=t_0$. Observe that there is no change in stability or basin size of the QSE as $t$ changes. Fig.~\ref{fig:hopf_shift_rate} shows typical trajectories starting at an arbitrary initial condition within the basin of attraction using fixed $\Delta$ and two values of $\rho$. As shown in the diagram, there is a critical value $\rho_c$ such that for $\rho<\rho_c$ the system can track the QSE while for $\rho>\rho_c$ a tipping occurs near $t=t_0$.

Jan Sieber ({\em pers. comm.}) has pointed out that this case may still be quantifiable by numerical approximation of the $\rho_c$ that gives a heteroclinic connection from an (initial) saddle equilibrium at $(z,\lambda)=(0,0)$ to  a saddle periodic orbit at $(|z-\Delta|,\lambda)=(1,\Delta)$ for  the extended system (\ref{eq:hopf}) and (\ref{eq:shift}). Such a connection indicates  $\rho_c$ for which the (initial) saddle equilibrium moves away from the basin boundary of the stable equilibrium at $(z,\lambda)=(\Delta,\Delta)$.

\subsection{A fast-slow system with $R$-tipping}

A particularly interesting case of R-tipping can occur
in slow-fast systems that have a (unique, globally stable) QSE near locally folded critical (slow) manifold, of which the recently studied compost-bomb instability is a representative~\cite{Luke10,Wieczorek_etal_2010}.
Here, we consider a simple example
\begin{eqnarray}
\label{eq:sf1}
\epsilon\,\frac{dx_1}{dt}&=&x_2+\lambda+x_1(x_1-1),\\
\label{eq:sf2}
\frac{dx_2}{dt}&=&-\sum_{n=1}^N x_1^n,
\end{eqnarray}
with odd $N$, fast variable $x_1\in\mathbb{R}$, slow variable $x_2\in\mathbb{R}$, 
and small parameter $0<\epsilon\ll 1$. A unique equilibrium for~(\ref{eq:sf1})--(\ref{eq:sf2}), 
$
\tilde{x}(\lambda)=(0,-\lambda),
$
is asymptotically stable for any fixed value of $\lambda$, and globally asymptotically stable if $N\ge 5$.  
The slow dynamics is approximated by the one-dimensional {\em critical (slow) manifold},
$
S(\lambda)=\{(x_1,x_2)\in \mathbb{R}^{2}:x_2=-\lambda - x_1(x_1-1)\},
$
that has a fold,
$
L(\lambda)=(\frac{1}{2},-\lambda+\frac{1}{4}),
$
tangent to the fast $x_1$ direction. If $N\ge 5$, the fold defines a tipping threshold that is not associated with any basin boundary. Here, $S(\lambda)$ is partitioned into the attracting part, $S_a(\lambda$) for $x_1<\frac{1}{2}$, fold $L(\lambda)$ for $x_1=\frac{1}{2}$, and  repelling part, $S_r(\lambda)$ for $x_1>\frac{1}{2}$.

\subsubsection{The slow-fast system with steady drift}

\begin{figure}[t!]
  \begin{center}
    \includegraphics[type=eps,ext=.eps,read=.eps,width=13cm]
    {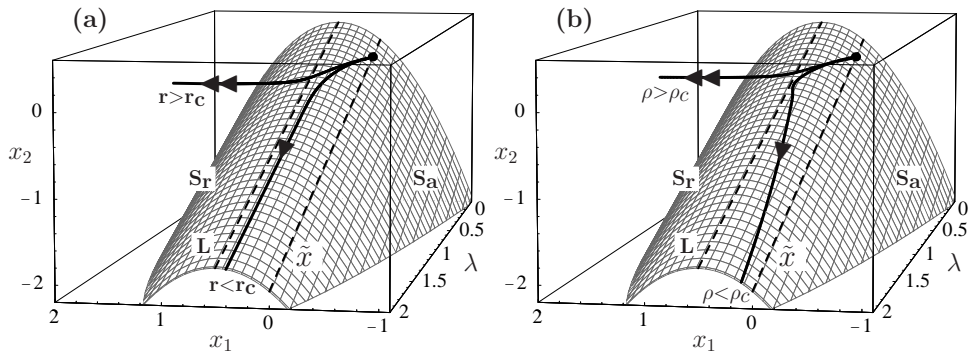}
    \caption{\label{fig:sw}
R-tipping in slow-fast systems with a unique quasi-stable equilibrium, $\tilde{x}$, and (gray surface) locally folded critical (slow) manifold, $S=S_a\cup L\cup S_r$, for (a) the steady drift problem~(\ref{eq:sf1})--(\ref{eq:rate}) and (b) the unsteady drift problem~(\ref{eq:sf1})--(\ref{eq:sf2}) and  (\ref{eq:nrate}), where $\epsilon=0.01$ and $N=1$. In (a), Eq.~(\ref{eq:rcsfs}) gives  $r_c=\frac{1}{2}$ and shown are trajectories for $r=0.4<r_c$ and $r=0.6>r_c$. In (b), Eq.~(\ref{eq:rcsfu}) gives $\rho_c\approx 0.99$ for the initial condition at the origin (black dot) and shown are trajectories for $\rho=0.7<\rho_c$ and $\rho=1>\rho_c$.
}
  \end{center}
\end{figure}

Consider~(\ref{eq:sf1})--(\ref{eq:sf2}) with a uniform drift of the QSE, $\tilde{x}(\lambda(t))$,  in the negative $x_2$ direction at a constant rate
\begin{equation}
\label{eq:rate}
\frac{d\lambda}{dt}=r>0,
\end{equation}
so that $\lambda$ becomes the second slow variable.
There is a critical rate, $r_c$, at which~(\ref{eq:sf1})--(\ref{eq:rate}) is destabilized, meaning that trajectories diverge away from the QSE for $r>r_c$. 
We can find this critical rate in the singular limit, $\epsilon\rightarrow 0$, by setting $\epsilon=0$ in~(\ref{eq:sf1}), differentiating the resulting algebraic equation with respect to $t$, and studying the {\em  projected reduced system}~\cite{FEN79}: 
\begin{eqnarray}
\label{eq:srs1}
\frac{dx_1}{dt}&=&\left(-r+\sum_{n=1}^N x_1^n\right)(2x_1-1)^{-1},\\
\label{eq:srs2}
\frac{d\lambda}{dt}&=&r,
\end{eqnarray}
that aproximates the slow dynamics for~(\ref{eq:sf1})--(\ref{eq:rate}) on the two-dimensional critical manifold,
$
S=\{(x_1,x_2,\lambda)\in \mathbb{R}^{3}:x_2=-\lambda - x_1(x_1-1)\}
$ 
(gray surface in  Fig.~\ref{fig:sw}).
Although~(\ref{eq:srs1})--(\ref{eq:srs2}) is typically singular at the one-dimensional fold,
$
L=\{(x_1,x_2,\lambda)\in \mathbb{R}^{3}:x_1=\frac{1}{2},\;\;x_2=-\lambda+\frac{1}{4}\},
$ 
its phase portrait can be constructed by rescaling time
$$
\frac{dt}{d\tau}=-(2x_1-1)\;\;\;\;\;\Rightarrow\;\;\;\;\; t=-\int_0^{\tau}\left(2x_1(s)-1\right)ds,
$$
producing the phase portrait for the {\em desingularised system}~\cite{KRU01}:
\begin{eqnarray}
\label{eq:ds1}
\frac{dx_1}{d\tau}&=& r -\sum_{n=1}^N x_1^n,\\
\label{eq:ds2}
\frac{d\lambda}{d\tau}&=&-r(2x_1-1),
\end{eqnarray}
and then reversing the direction of time on the repelling part of the critical manifold, $S_r$.
In this way, we find that for $0<r<\sum_{n=1}^N 2^{-n}$ trajectories for all initial conditions 
within $S_a$ converge to a stable invariant line that is defined by a constant $x_1$ satisfying $r=\sum_{n=1}^N x_1^{n}$, meaning that trajectories remain close to the QSE, $\tilde{x}$, for all time [$r<r_c$ in Fig.~\ref{fig:sw}(a) ]. However,  for $r>\sum_{n=1}^N 2^{-n}$, trajectories for all initial conditions within $S_a$ reach the fold, $L$, where they ``slip off'' the critical manifold and diverge away from the QSE in the fast $x_1$ direction [$r>r_c$ in Fig.~\ref{fig:sw}(a)]. Hence, system~(\ref{eq:sf1})--(\ref{eq:rate}) exhibits R-tipping and, for $\epsilon\rightarrow 0$, the critical rate is
\begin{equation}
\label{eq:rcsfs}
r_c=\sum_{n=1}^N 2^{-n}.
\end{equation}

\subsubsection{The slow-fast system with unsteady drift}

We now consider~(\ref{eq:sf1})--(\ref{eq:sf2}) with a nonuniform drift
\begin{equation}
\label{eq:nrate}
\frac{d\lambda}{dt}=\rho \,e^{-\lambda}
\end{equation}
that is a logarithmic growth, $\lambda(t)=\ln\left[\rho(t-t_0)+e^{\lambda(t_0)}\right]$, where we assume $\rho>0$. Again, there is a critical rate, $\rho_c$, at which system~(\ref{eq:sf1})--(\ref{eq:sf2}) and (\ref{eq:nrate}) is destabilized. The key difference from the steady drift problem is that $\rho_c$ depends on the initial condition within $S_a$. This is because the desingularised system
\begin{eqnarray}
\label{eq:ds11}
\frac{dx_1}{d\tau}&=& e^{-\lambda}\,\rho -\sum_{n=1}^N x_1^n,\\
\label{eq:ds22}
\frac{d\lambda}{d\tau}&=&-e^{-\lambda}\,\rho\,(2x_1-1),
\end{eqnarray}
has a saddle equilibrium for all $\rho>0$, corresponding to a {\em folded-saddle singularity}~\cite{ARN94,SZM01}:
$$
F=\left(x_{1,F},\lambda_F(\rho)\right)=\left(\frac{1}{2},-\ln\sum_{n=1}^N\,2^{-n}/\rho\right),
$$
for the projected reduced system.
One can use the theory developed in~\cite[Sec.4]{Wieczorek_etal_2010} to approximate the critical value, $\rho_c$. Given $F$, the eigenvector
$$
w=\left(-q/p+\sqrt{2+(q/p)^2},1\right)^T
$$
corresponding to the stable eigendirection of the saddle $F$ for~(\ref{eq:ds11})--(\ref{eq:ds22}), an initial condition $(x_{1,0},\lambda_0)$  within $S_a$, and as far as $\epsilon\rightarrow 0$, the critical rate can be calculated using~\cite[Eq.(4.12)]{Wieczorek_etal_2010}
to give
\begin{equation}
\label{eq:rcsfu}
\rho_c\approx p\,\,\exp{\left(
\lambda_0+\frac{1/2-x_{1,0}}{-q/p+\sqrt{2+(q/p)^2}}
\right)},
\end{equation}
where $p=\sum_{n=1}^N\,2^{-n}$, $q=\sum_{n=1}^N\,n2^{-n}$.
Below the critical rate, the trajectory misses the fold, $L$, and approaches the QSE, $\tilde{x}$, as time tends to infinity [$\rho<\rho_c$ in Fig.~\ref{fig:sw}(b)]. Above the critical rate, the trajectory reaches $L$ and diverges from the QSE in the fast $x_1$ direction [$\rho>\rho_c$ Fig.~\ref{fig:sw}(b)].
Note that in this example, the critical rate of parameter variation is of the same order as the {\em slow} dynamics - only when the parameter variation is {\em very slow} with respect to the slow variable and there are three timescales is tracking guaranteed. In this sense, the rate dependent tipping occurs when the slow and very slow timescales are no longer separable.

\section{B-, N- and R-tipping examples in a simple climate model}
\label{sec:BNR-tipping_egs}

We present a simple climate model that independently show, under differing circumstances, all three types of tipping. In its deterministic version, this is a
``zero dimensional'' global energy balance model originally introduced by
Fraedrich~\cite{Fraedrich1979}:
\begin{equation}
  \label{eq:imbalance}
  c\frac{dT}{dt}=R\downarrow-R\uparrow.
\end{equation}
The state variable $T$ represents an average surface temperature of an ocean
on a spherical planet subject to radiative heating.  Eq.~\eqref{eq:imbalance}
is a deterministic energy conservation law where the constant $c$ represents
the thermal capacity of a well-mixed ocean layer of depth $30$m covering
70.8\% of the Earth's surface. The incoming solar radiation $R\downarrow$ and
outgoing radiation $R\uparrow$ are modelled as
\begin{equation*}
  \label{eq:incoming}
  R\downarrow=\frac{1}{4}\mu I_0(1-\alpha_p(T)),
  \quad
  R\uparrow=e_{SA}\sigma T^4.
\end{equation*}
Here $I_0$ is the solar constant and the parameter $\mu$ allows for
variations in the planetary orbit, or in the solar constant. An ice-albedo
feedback is introduced to link variations in temperature with changes of ice
and thus of albedo $\alpha_p$: Fraedrich~\cite{Fraedrich1979} uses a quadratic
relation
\begin{equation}
  \label{eq:albedo}
  \alpha_p(T)=a_2-b_2T^2,
\end{equation}
where the parameters $a_2>1$ and $b_2$ control the albedo magnitude and slope
of the albedo-temperature relation.  The outgoing radiation term is obtained
by the Stefan-Boltzmann law, where $e_{SA}$ is the effective emissivity and
$\sigma$ is the Stefan-Boltzmann constant. With these choices
\eqref{eq:imbalance} is written as~\cite[Eq. 4.1]{Fraedrich1979}:
\begin{gather}
  \label{eq:Fraedrich}
  \frac{dT}{dt}=f(T)=c^{-1}a(-T^4+b_\mu T^2-d_\mu),\\
  a=e_{SA}\sigma/c,\quad
  b_\mu=\mu I_0b_2/4e_{SA}\sigma,\quad
  d_\mu=-\mu I_0(1-a_2)/4e_{SA}\sigma.\nonumber
\end{gather}
Table~\ref{tab:fraed_values} shows the values of constants and parameters for
the system at equilibrium.
\begin{table}
  \centering
  \begin{tabular}{ll|ll}
  \hline
    $I_0$ & $1366$ W\,m$^{-2}$ & $\mu$ & 1  \\
    $\sigma$ & $5.6704\times 10^{-8}$ W\,m$^{-2}$\,K$^{-4}$ & $b_2$ & $1.690\times 10^{-5}$ K$^{-2}$\\
    $c$ & $10^{8}$ kg\,K\,s$^{-2}$ & $a_2$ & $1.6927$\\
    $e_{SA}$ & 0.62  & & \\
    \hline
  \end{tabular}
  \caption{Values of the constants and parameters for Eq.~\eqref{eq:Fraedrich}.}
  \label{tab:fraed_values}
\end{table}
Sutera~\cite{Sutera1981} reformulates Fraedrich's model to incorporate
stochastic forcing:
\begin{equation}
  \label{eq:Sutera}
  dT=f(T)dt+\sqrt{\nu}\,dW,
\end{equation}
with $f(T)$ as in~\eqref{eq:Fraedrich}, where $dW$ is a normalised Wiener
(white noise) process such that $(dW)^2$ has dimension of time $t$, $\nu$ has
dimension of $1/t$ and the variance of $\sqrt{\nu}\,dW$ per unit time is
$\nu$.

For $\mu$ larger than a critical value $0<\mu_c<1$, the deterministic
system~\eqref{eq:Fraedrich} has two equilibria $T^+$ (stable) and $T^-$
(unstable). A saddle-node bifurcation takes place at some $\mu=\mu_c$ with $0<\mu_c<1$, where the
two equilibria $T^\pm$ merge and disappear. Sutera~\cite{Sutera1981} studies
N-tipping in the stochastically forced system~\eqref{eq:Sutera} for
$\mu>\mu_c$, as a function of the distance $\mu-\mu_c$ from the bifurcation
value. Namely, they compute the exit time such that the process jumps over the
``potential barrier'' $T^-$ and falls irreversibly to `ice-covered earth'.

\begin{figure}%
\begin{center}
\includegraphics[width=12cm,clip=]{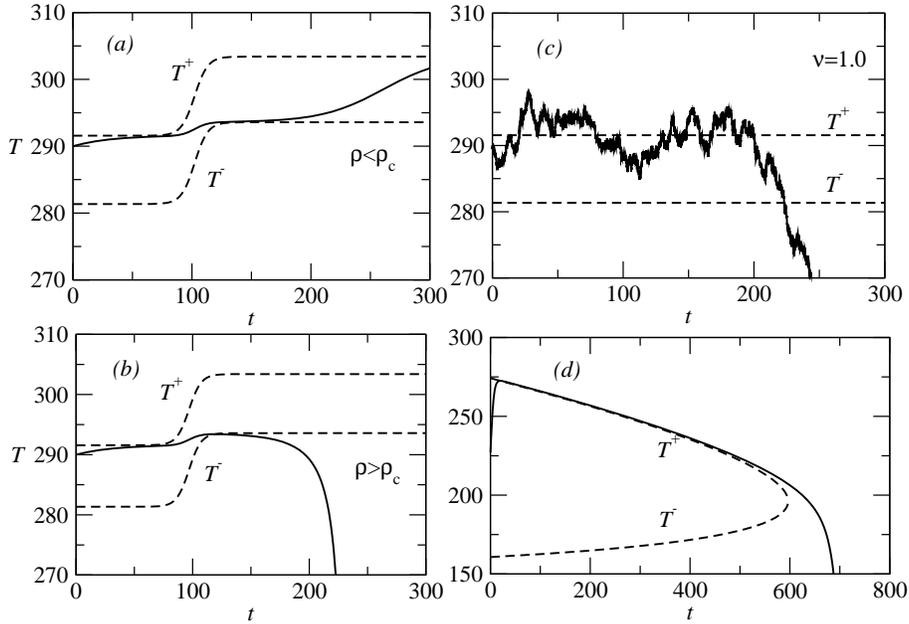}%
\end{center}
\caption{ Illustrations of trajectories for the Sutera-Fraedrich model
  (\ref{eq:Sutera}) showing the presence of all three tipping types for
  parameters in Table~\ref{tab:fraed_tipping} - horizontal axis (years) 
  vertical axis (Kelvin). The solid lines show 
  system trajectories while the dash lines show the location of the QSE -- the 
  branch $T^+$ is stable while the branch $T_-$ is unstable in this 
  model. Panels (a)-(b) show
  an R-tipping for a smooth change of parameters between two steady states. In 
  (a) $\rho=0.18$ the system 
  returns to the QSE after a transient. In (b) $\rho=0.19$ the system becomes 
  unbounded indicating a critical value $\rho_c\approx 0.185$ yr$^{-1}$. (c) 
  shows an example of N-tipping in the presence of noise
  of amplitude $\nu=1.0$ yr$^{-1}$; (d) shows an example of B-tipping on 
  decreasing $\mu$ uniformly from $1$ at a constant rate. Note that in case (d) 
  the two QSE coalesce at a saddle-node bifurcation.}%
\label{fig:fraed_tipping}%
\end{figure}

We illustrate in Fig.~\ref{fig:fraed_tipping} three situations where the
Sutera-Fraedrich model exhibits ``pure'' B-, N- and R-tipping independently; parameter
values are detailed in Table~\ref{tab:fraed_tipping}. Panels (a-b) show an example of R-tipping, (c) of N-tipping and (d) of B-tipping. For case (a-b), we evolve the dimensionless parameter $\lambda$ according to the ODE $d\lambda/dt=\rho\lambda(1-\lambda)$ and set $b_2=(1-\lambda)b_2^{\mathrm{init}}+\lambda b_2^{\mathrm{final}}$ - for this figure we use initial values $\lambda=10^{6}$ and $T=290$K. The value of $a_2$ is calculated to ensure that $b_\mu^2-4d_\mu$ is held constant for the parameter groups defined in (\ref{eq:Fraedrich}). The constant $\rho$ can be thought of as simply scaling the rate of passage from the initial to the final values given in the Table.

\begin{table}%
\begin{center}
\begin{tabular}{lccc}
Parameter & (a-b) & (c) & (d)\\
\hline
$\mu$  & 1.0 & 1.0 & decreases from $1.0$ \\
& & & at rate $-0.0004$ yr$^{-1}$ \\
$b_2$ (K${}^{-2}$) & initial $1.690\times 10^{-5}$  & $1.690\times 10^{-5}$ & $1.04\times 10^{-5}$ \\ 
& final $1.8350\times 10^{-5}$ & & \\
$a_2$ & initial $ 1.6927$  &  $ 1.6927$ & $1.2$\\
& final $1.8168$ & & \\
$\nu$ (yr$^{-1}$)& 0 & 1.0 & 0 \\
\hline
\end{tabular}
\end{center}
\caption{ Parameter values for simulations shown in Fig.~\ref{fig:fraed_tipping}. For case (a-b) we interpolate between the values given along the curve such that $b_{\mu}^2-4d_{\mu}$ is constant at a rate proportional to $\rho$. In case (c) all parameters are fixed but noise is added, while in case (d) we impose a steady drift of the parameter $\mu$ downwards.}%
\label{tab:fraed_tipping}%
\end{table}

\section{Summary and conclusions}
\label{sec:discuss}

It is of great practical importance to understand the theoretical mechanisms behind tipping phenomena in the climate system as well as other
systems. We have proposed here that such mechanisms can be effectively divided
into three distinct categories: bifurcation-induced, noise-induced and rate-dependent tipping,
respectively denoted as B-, N- and R-tipping. In particular, we describe R-tipping, a mechanism that may be exhibited by
(sub-systems of) the climate system independently of the presence or absence
of the other types of tipping.

\begin{figure}%
\begin{center}
\includegraphics[width=12cm,clip=]{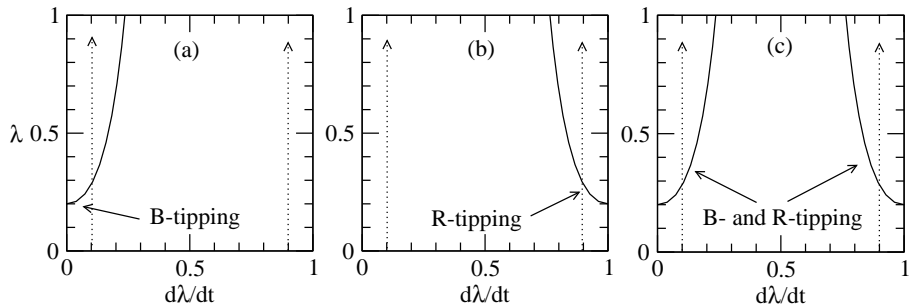}%
\end{center}
\caption{Schematic illustration of different possible system behaviours on
  ramping the parameter $\lambda$ at differing rates (dashed arrows) through
  the region $\lambda\in[0,1]$. (a) shows an example where there is a
  B-tipping for low rates (quasi-static) that disappears for high enough
  rates. (b) shows an example where there is no tipping for small enough rates
  but R-tipping for large enough rates. (c) shows both B- and R-tipping, but
  there is a range of rates where no instability appears.}%
\label{fig:final}%
\end{figure}

In realistic models, tipping effects may be associated with a combination of
the three mechanisms, and it will be a challenge to understand
this more general case. For example B-tipping, usually associated with slow
changes in a parameter, may turn into R-tipping upon increasing the rate of
change for the parameters. However, as schematically illustrated in
Fig.~\ref{fig:final}, completely new mechanisms may appear on increasing the
rate, including the possibility that B-tipping may be suppressed for fast enough variation of parameters. Alternatively, the B-tipping may persist but an  R-tipping mechanism may come into play before the B-tipping is reached. 

We emphasise that neither N-tipping and R-tipping require any change of
stability. Hence there is no reason to assume that the techniques
of~\cite{Dakos_etal_2008}, based on a de-trended autoregressive model for B-tipping, should deliver useful predictions in such cases - as noted
by~\cite{Ditlevsen_etal_2010}, N-tipping is intrinsically unpredictable. We
are investigating whether any novel predictive technique may be
developed for R-tipping. Those cases of R-tipping that can be reduced to a local bifurcation 
in a co-moving system may be expected to be predictable using similar methods; this includes 
the examples in Sec.~3(a),(b)(i) with steady drift. 
In more complex cases, $\rho_c$ may still be quantifiable by global heteroclinic 
bifurcations for an extended system, for example (\ref{eq:nf2}) and (\ref{eq:shift}) or (\ref{eq:hopf}) and (\ref{eq:shift}) in Sec.~3(b)(ii).

The classification proposed here may be applicable to a wide range of
open systems under the influence of noise and/or parameter changes. Recent work
of Nene and Zaikin \cite{Nene_Zaikin_2010} suggests there may be interesting
applications of rate-dependent bifurcation theory to determine cell
fate. There are potentially many other application areas, from mechanics and ecology to economics and social sciences where tipping points are of interest. We suggest that this will be an area of significant mathematical development in the coming years.

~

\noindent {\bf Acknowledgements:}
  The Authors thank the Isaac Newton Institute for hosting the programme
  ``Mathematical and Statistical Approaches to Climate Modelling and
  Prediction'' in autumn 2010, where this topic was initially discussed, and
  are indebted to Alexei Zaikin and Jan Sieber for stimulating conversations
  in relation to this work.

\bibliographystyle{plain}

\end{document}